\documentclass[12pt]{article}
\usepackage{latexsym,amsmath,amsopn,amssymb,amsthm,amsfonts,epsf}

\begin{document}

\vspace*{3cm} {\normalsize\bf About construction of orthogonal
wavelets with compact support and with scaling coefficient  N }
\vspace{28pt}

{\normalsize P.~N.~Podkur and N.~K.~Smolentsev }

\vspace{28pt}

\thispagestyle{empty} \baselineskip=13pt
\begin{quote}

In paper \cite {BEJ} with using of the Cuntz algebra
representation some methods of construction of wavelets with
scaling coefficient $N\ge 2$ are considered. In paper \cite {BJ}
it is shown, a construction of wavelets at the prescribed scaling
function $\varphi (x)$. In this paper a simple method of
construction of scaling function $\varphi (x)$ and orthogonal
wavelets with the compact support for any natural coefficient of
scaling $N\ge 2$ is given. Examples of construction of wavelets
for coefficients of scaling $N=2$ and $N=3$ are produced.
\end{quote}

{\bf 1. Scaling functions and wavelets.} Let $N\ge 2$ is an
integer, $\mathbb {Z}$ is set of all integers and $L^2 (\mathbb
{R})$ is Hilbert space of square integrable functions.

{\bf Definition 1.} {\it Function  $\varphi(x) \in
L^2(\mathbb{R})$ is called $N$-scaling, if it can be represented as
$$
\varphi(x)=\sqrt{N}\sum_{n\in \mathbb{Z}} h_n \varphi(Nx-n),
\eqno(1)
$$
where coefficients $h_n$,\ $n\in \mathbb{Z}$ satisfy to condition
$\sum_{n}|h_n|^2<\infty$. The relationship (1) is called the
$N$-scale equation (refinement equation). The set $\{h_n\}$ of
coefficients of expansion in the equation (1) is called the
scaling filter.}

{\bf Note 1.} If $N$-scaling function $\varphi(x)$ has the compact
support of length $L$, then the sum in equation (1) is finite,
contained at most $L(N-1)+1$ components.

The Fourier transform of $N$-scale equation is
$$
\widehat{\varphi}(\omega)=H_0\left( \frac{\omega}{N}\right)
\widehat{\varphi}\left(\frac{\omega}{N} \right) , \eqno(2)
$$
where
$$
H_0(\omega)=\frac {1}{\sqrt{N}} \sum_{n\in \mathbb{Z}} h_n
e^{-in\omega}. \eqno(3)
$$
The function $H_0(\omega)$ is called {\it frequency function} of
scaling function $\varphi(x)$.

In the orthogonal case translations of scaling function $\varphi
(x-n), \ n\in \mathbb {Z}$ form orthonormal basis of the subspace
$V_0$ in $L^2 (\mathbb {R})$, and translations $ \varphi _ {1, n}
(x) = \sqrt {N} \varphi (Nx-n), \ n\in \mathbb {Z}$ on $1/N $,
form orthonormal basis of the subspace $V_1 $ in $L^2 (\mathbb
{R})$. Thus $V_0 \subset V_1$. In the orthogonal case to the
scaling function $\varphi (x)$ corresponds $N-1$
wavelets-functions $\psi^1 (x)$..., $\psi ^ {N-1} (x)$, for each
of which translations  $\psi^k _ {0, n} (x) = \psi^k (x-n), \ n\in
\mathbb {Z}$ form orthonormal basis of subspaces $W_0^k$ in $L^2
(\mathbb {R})$, and expansion in the direct sum of orthogonal
subspaces $V_1=V_0\oplus W_0^1\oplus \dots \oplus W_0 ^ {N-1}$ be
valid.

Wavelets $\psi^1 (x) $..., $ \psi ^ {N-1} (x)$ form orthonormal
basis $L^2 (\mathbb {R})$:
$$
\{\psi^k_{j,n}(x)=\sqrt{N^j}\psi^k(N^jx-n),\ j,n\in \mathbb{Z},\
k=1,2,\dots , N-1 \}.
$$
As wavelets $\psi^1 (x)$..., $\psi ^ {N-1} (x)$
belong to space $V_1 $ they are decomposed on basis of this space,
$$
\psi^k(x)=\sqrt{N}\sum_{n\in \mathbb{Z}} g^k_n \varphi(Nx-n).
\eqno(4)
$$
The coefficients $\{g^k_n \}$ is called {\it filters of wavelets}
$\psi^k (x)$, \ $k=1,2, \dots, N-1$. Let
$$
H_k(\omega)=\frac {1}{\sqrt{N}} \sum_{n\in \mathbb{Z}} g^k_n
e^{-in\omega} \eqno(5)
$$
-- the frequency functions corresponding to wavelets $\psi^1 (x)
$..., $ \psi ^ {N-1} (x) $. The Fourier transform of equalities
(4) is
$$
\widehat{\psi}^k(\omega)=H_k\left( \frac{\omega}{N}\right)
\widehat{\varphi}\left(\frac{\omega}{N} \right) .
$$

For the frequency functions $H_k (\omega)$ the following matrix is
unitary  \cite {BJ}, \cite {Db},
$$
H(z)=\left(
\begin{array}{cccc}
  H_0(z) & H_0(\rho z) & \ldots & H_0(\rho^{N-1} z) \\
  H_1(z) & H_1(\rho z) & \ldots & H_1(\rho^{N-1} z) \\
  \hdotsfor[2.5]{4} \\
  H_{N-1}(z) & H_{N-1}(\rho z) & \ldots & H_{N-1}(\rho^{N-1} z) \\
\end{array} \right),  \eqno(6)
$$
where $z=e^{-i\omega}$ and $\rho=e^{-i2\pi/N}$.
The matrix (6) has special view. It is possible to avoid of
this special view of the matrix $H(z)$ with Fourier transform on
cyclic group $\mathbb{Z}/N\mathbb{Z} = \{1, \rho, \rho^2,\dots, \rho^{N-1}\}$
\cite{BJ}. We shall define
$$
A_{k,j}(w)=\frac {1}{\sqrt{N}}\sum_{z^N=w}z^{-j}H_k(z). \eqno(7)
$$
It is easy to see, that the sum on the right depends from $w = z^N$.
Also transformation (7) accurate within coefficient $\sqrt {N}$
is sample of elements with degrees $z ^ {kN}$ in polynomials
$H_k (z)$, $z ^ {-1} H_k (z) $..., $z ^ {-N+1} H_k (z)$.
Inverse transformation is defined by the formula \cite {BJ}
$$
H_k(z)=\frac {1}{\sqrt{N}}\sum_{j=0}^{N-1}z^{j}A_{k,j}(z^N).
\eqno(8)
$$
From last relation we shall obtained the following matrix equality:
$$
H(z)=\frac {1}{\sqrt{N}}A(z^N)\left(
\begin{array}{cccc}
  1 & 1 & \ldots & 1 \\
  z & \rho z & \ldots & \rho^{N-1} z \\
  \hdotsfor[2.5]{4} \\
  z^{N-1} & \rho^{N-1} z^{N-1} & \ldots & \rho^{((N-1)^2)} z^{N-1} \\
\end{array}
\right) = A(z^N)R(z). \eqno(9)
$$
In this expression the matrix $A (z^N)$ is already arbitrary
unitary matrix with polynomial elements. Now specificity of the matrix $H(z)$
go to the matrix
$$
R(z)=\frac {1}{\sqrt{N}}\left(
\begin{array}{cccc}
  1 & 1 & \ldots & 1 \\
  z & \rho z & \ldots & \rho^{N-1} z \\
  \hdotsfor[2.5]{4} \\
  z^{N-1} & \rho^{N-1} z^{N-1} & \ldots & \rho^{(N-1)^2} z^{N-1} \\
\end{array}
\right). \eqno(10)
$$
Let's mark, that the matrix $R(z)$ is unitary on the unit circle
$z=e^{-i\omega}$.

Specifying the polyphase matrix $A(w)$, we can construct the
matrix of frequency functions $H(z)$ by the formula (8) and,
together with it, frequency functions of wavelets $H_1(z), \dots, H_{N-1}(z)$,
hence, and wavelets $\psi^1(x)\dots , \psi^{N-1}(x)$.

In work \cite {BJ} the scheme of construction of the polyphase
matrix $A (z^N)$ is given in the supposition, that polynomial
frequency function $H_0 (z)$ is prescribed. Then it is possible to consider,
that the first row of the matrix $A_{0j}(z^N) $ is known,
$$
A_{0,j}(w)=\frac {1}{\sqrt{N}}\sum_{z^N=w}z^{-j}H_0(z), \eqno(11)
$$
and it is necessary to construct remaining row of the matrix $A(w)$.

In the given work we shall give the simple scheme of construction
of the unitary matrix $A(w)$ which elements are polynomials with
real coefficients. It allows to define both the scaling function
$\varphi (x)$ with compact support and with scaling coefficient
$N>2$, and orthogonal wavelets $\psi^1(x)\dots , \psi^{N-1}(x)$.


{\bf 2. Scheme of wavelets construction.} From above constructions
and methods of work \cite {BJ} follows that orthogonal systems of
wavelets can be determine by the unitary matrix $A(w)$ with
polynomial elements with using of the formula $H(z) = A(z^N)
R(z)$, where $R(z) $ -- the special matrix (10). We shall give the
simple method of construction enough big set of unitary matrixes
$A(w)$ with polynomial elements. It will allow to obtain both the
$N$-scaling function with the compact support, and orthogonal
wavelets.

Let's choose any orthogonal matrix
$A_0 = \{a _ {ij}, \ i, j = 0,1,\dots, N-1 \}$ of the order
$N\ge 2 $. We shall multiply it on the diagonal unitary matrix
$D_k (w) = {\rm diag} (w ^ {k_0}, w ^ {k_1}, \dots, w ^ {k _ {N-1}})$,
where $k = (k_0, k_1, \dots, k_{N-1})$ is set of integers and
$|w | = 1$, and then -- on the orthogonal matrix
$B_0 = \{b _ {ij}, \ i, j = 0,1, \dots, N-1 \}$.
In outcome we shall obtain unitary matrix
$$
A(w) = A_0D_k(w)B_0, \eqno(12)
$$
which elements, $A_{ij}=\sum_{s=0}^{N-1}a_{is}b_{sj}w^{k_s}$,
are polynomials on the variable $w$ with real coefficients.

Now we shall substitute $w = z^N$, where $z=e^{-i\omega}$. We
shall obtain the unitary matrix $A (z^N)$ with polynomial elements
and real coefficients. We shall multiply it on the unitary matrix
$R(z)$. Then we shall obtain the unitary matrix $H (z) $ of
frequency polynomial functions $H_0(z)$, $H_1(z) $...,
$H_{N-1}(z)$ of orthogonal system of wavelets $ \varphi (x)$, $
\psi^1(x) \dots, \psi^{N-1}(x)$, where the first function $\varphi
(x)$ is scaling, and remaining -- wavelets. Thus,
$$
H(z)=\left(
\begin{array}{cccc}
  H_0(z) & H_0(\rho z) & \ldots & H_0(\rho^{N-1} z) \\
  H_1(z) & H_1(\rho z) & \ldots & H_1(\rho^{N-1} z) \\
  \hdotsfor[2.5]{4} \\
  H_{N-1}(z) & H_{N-1}(\rho z) & \ldots & H_{N-1}(\rho^{N-1} z) \\
\end{array} \right)=A_0D_k(z^N)B_0R(z).  \eqno(13)
$$
From (13) follows the expression for frequency functions:
$$
H_{k}(z)=\frac{1}{\sqrt{N}}\sum_{s,j=0}^{N-1}a_{ks}b_{sj}z^{j}z^{Nk_s},\
k=0,1,\dots, N-1. \eqno(14)
$$

In order to the obtained the functions $H_k(z)$ would be frequency
functions of orthogonal wavelets, it is necessary, that the sum of
coefficients for $H_0 (z)$ would be equal to unit, and the sums of
coefficients for remaining functions $H_1(z)$..., $H_{N-1}(z)$
would be equal to zero:
$$
\frac{1}{\sqrt{N}}\sum_{s,j=0}^{N-1}a_{0s}b_{sj}=
\frac{1}{\sqrt{N}}\sum_{s=0}^{N-1}a_{0s}\sum_{j=0}^{N-1}b_{sj} =1,
$$
$$
\frac{1}{\sqrt{N}}\sum_{s,j=0}^{N-1}a_{ks}b_{sj}=
\frac{1}{\sqrt{N}}\sum_{s=0}^{N-1}a_{ks}\sum_{j=0}^{N-1}b_{sj} =0,
\qquad k=0,1,\dots, N-1.
$$

These equalities can be represented in the matrix view:
$$
\left(
\begin{array}{cccc}
  a_{00} & a_{01} & \ldots & a_{0,N-1} \\
  a_{10} & a_{11} & \ldots & a_{1,N-1} \\
  \hdotsfor[2.5]{4} \\
  a_{N-1,0} & a_{N-1,1} & \ldots & a_{N-1,,N-1} \\
\end{array} \right)
\left(
\begin{array}{c}
  b_{00}+  \dots + b_{0,N-1} \\
  b_{10}+  \dots + b_{1,N-1} \\
  \ldots \ldots \ldots\\
  b_{N-1,0}+  \dots + b_{N-1,,N-1} \\
\end{array} \right)=
\left(
\begin{array}{c}
  \sqrt{N} \\
  0 \\
  \ldots \\
  0 \\
\end{array} \right).  \eqno(15)
$$

Choosing various orthogonal matrixes $A_0$ and $B_0$,
which satisfy the equality (15), we obtain various frequency
functions of wavelets (14).

For construction enough simple class of orthogonal wavelets
with the compact support and scaling coefficient $N>2$,
we shall take as an orthogonal matrix $A_0$ the following matrix:
$$
A_0=\left(
\begin{array}{ccccc}
  1/\sqrt{N} & 1/\sqrt{N} & 1/\sqrt{N} &  \ldots & 1/\sqrt{N} \\
  1/\sqrt{2} & -1/\sqrt{2} & 0 & \ldots & 0  \\
  1/\sqrt{6} & 1/\sqrt{6} & -2/\sqrt{6} & \ldots & 0  \\
  \hdotsfor[2.5]{5} \\
   1/\sqrt{N(N-1)} & 1/\sqrt{N(N-1)} & 1/\sqrt{N(N-1)} &  \ldots & -(N-1)/\sqrt{N(N-1)} \\
\end{array} \right).
$$

The matrix $A_0$ transform vector of units $e = (1,1, \dots, 1) $
to the vector $\sqrt {N} e_0 = (\sqrt {N}, 0, \dots, 0)$, \ $A_0 e
= \sqrt {N} e_0 $. Then from equality (15) follows, that elements
of the orthogonal matrix $B_0$ should satisfy to the following
system of equations:
$$
\left\{
\begin{array}{ccc}
  b_{00}+b_{01}+ \ldots + b_{0,N-1} & = & 1 \\
  b_{10}+b_{11}+ \ldots + b_{1,N-1} & = & 1  \\
  \hdotsfor[2.5]{3} \\
  b_{N-1,0}+b_{N-1,1}+ \ldots + b_{N-1,N-1} & = & 1  \\
\end{array} \right. . \eqno(16)
$$

The solution of this system will be any set of orthonormal vectors (rows)
which coordinates satisfy to the equation of the plane
$x_0 + x_1 + \ldots + x _ {N-1} =1 $ in $ \mathbb {R} ^N $. It is obvious,
that coordinates of basis vectors
$e_0 = (1, 0, \dots, 0) $, $e_1 = (0, 1, 0, \dots, 0)$...,$e_{N-1} = (0,
\dots, 0, 1)$ satisfy to this equation. The given solution corresponds to
the identity matrix $B_0$. Any other solution can be obtained by rotation of
the basis solution $e_0, e_1, \dots, e_{N-1}$ around of vector
$e = e_0+e_1 + \dots + e_{N-1} $, i.e. in the plane
$x_0 + x_1 + \ldots + x_{N-1} =1 $. We shall find these solutions.
We shall take rotation around of axis $Ox_0$:

$$
M=\left(
\begin{array}{ccccc}
  1 & 0 & 0 & \ldots & 0 \\
  0 & m_1^1 & m_2^1 & \ldots & m_{N-1}^1  \\
  \hdotsfor[2.5]{5} \\
  0 & m_1^{N-1} & m_2^{N-1} & \ldots & m_{N-1}^{N-1}  \\
\end{array} \right) . \eqno(17)
$$

As $A_0e = \sqrt {N} e_0 $, then rotation around of axis $e$ is
given by the matrix $M_e=A_0^{-1} MA_0$. Then rows of the matrix
$B_0$ will consist of coordinates of vectors-columns which are
obtained from $e_0, e_1, \dots, e_{N-1}$ by action on them matrix
$M_e$. Therefore the matrix $B_0$ is transposed to $M_e$. Then
$$
H_M (z) = A_0 D_k(z^N) M_e^T R(z)= A_0 D_k(z^N) A_0^{T}M^T A_0 R(z),
\eqno(18)
$$
where $M $ -- any orthogonal matrix of view (17) and $D_k(w)=
{\rm diag}(w^{k_0},w^{k_1},\dots, w^{k_{N-1}})$.

The formula (18) gives the direct method of construction the
big family of frequency functions
$H_0(z) $, $H_1(z) $..., $H_ {N-1}(z)$ and orthogonal wavelets
with the compact support
$\varphi(x)$, $ \psi^1(x) \dots, \psi^{N-1}(x)$. Wavelets of the
family depend of the orthogonal matrix $M $ of view (17)
and of the vector of degrees $k = (k_0, k_1, \dots..., k_{N-1})$
which it is possible to set arbitrarily.


{\bf 3. Construction of orthogonal wavelets with compact support
for $N=2$.} In the given section we shall show by the example of
scale $N=2$ effectiveness of the wavelets construction scheme
explained above. Though the matrix $D_k (w)$ can be anyone, we
shall take for example the diagonal matrix $D_1 (w)={\rm diag}
\{1, w \}$, $|w | = 1$. In case $N=2$ orthogonal matrixes $A_0 $
and $B_0 $ can be in the general view:
$$
A_0=\left(
\begin{array}{cc}
  \cos t & \sin t \\
  -\sin t & \cos t \\
\end{array}
\right), \qquad B_0=\left(
\begin{array}{cc}
  \cos u & \sin u \\
  -\sin u & \cos u \\
\end{array}
\right).
$$

Then
$$
H(z)=\frac{1}{\sqrt{2}}\left(\begin{array}{cc}
  \cos t & \sin t \\
  -\sin t & \cos t \\
\end{array}\right)
\left(\begin{array}{cc}
  1 & 0 \\
  0 & z^2 \\
\end{array}\right)
\left(\begin{array}{cc}
  \cos u & \sin u \\
  -\sin u & \cos u \\
\end{array}\right)
\left(\begin{array}{cc}
  1 & 1 \\
  z & \rho z \\
\end{array}\right).
$$

Frequency functions are
$$
H_0(z)=\frac{1}{\sqrt{2}}\left(\cos t \cos u +(\cos t \sin u)z
-(\sin t \sin u)z^2 + (\sin t \cos u)z^3\right), \eqno(19)
$$
$$
H_1(z)=\frac{1}{\sqrt{2}}\left(-\sin t \cos u -(\sin t \sin u)z
-(\cos t \sin u)z^2 + (\cos t \cos u)z^3\right), \eqno(20)
$$

The sum of coefficients of frequency function $H_0(z)$
should be equal to unit, and the sum of coefficients of frequency
function $H_1(z)$ should be equal to zero. The system (15) becomes:
$$
\left(\begin{array}{cc}
  \cos t & \sin t \\
  -\sin t & \cos t \\
\end{array}\right)
\left(\begin{array}{c}
  \cos u + \sin u \\
  \cos u - \sin u \\
\end{array}\right)=
\left(\begin{array}{c}
  \sqrt{2} \\
  0 \\
\end{array}\right),
$$
$$
\left\{\begin{array}{ccc}
  \cos u + \sin u & = & \sqrt{2}\cos t \\
  \cos u - \sin u & = & \sqrt{2}\sin t\\
\end{array}\right. .
$$
Solving last system, we obtain, $u=\pi/4-t$.

Thus, we have constructed the family of frequency functions of the
wavelets specified by formulas (19), (20) in which $u =\pi/4-t$.
After elimination of the variable $u$, we obtain::
$$
H_0(z)=\frac{1}{4}\left(1 +\cos 2t +\sin 2t +(1 +\cos 2t -\sin
2t)z + \qquad \qquad\qquad\qquad\right.
$$
$$
\qquad\qquad\qquad\qquad \left. +(1 -\cos 2t -\sin 2t)z^2 +(1
-\cos 2t +\sin 2t)z^3\right), \eqno(21)
$$
$$
H_1(z)=\frac{1}{4}\left(-1 +\cos 2t -\sin 2t +(1 -\cos 2t -\sin
2t)z + \qquad \qquad\qquad\qquad\right.
$$
$$
\qquad\qquad\qquad\qquad \left. +(-1 -\cos 2t +\sin 2t)z^2 +(1
+\cos 2t +\sin 2t)z^3\right).
 \eqno(22)
$$

The given frequency functions $H_0(z)$ and $H_1(z)$ coincide with
the same, but obtained other methods in work \cite {BEJ}. Various
wavelets of Haar, Daubechies wavelets and their analogs include
into this family. In the following section some examples are
given.

Choosing other matrix $D_k(z^N)$, similarly we can construct
other orthogonal wavelets with other support length.

{\bf 4. Examples of scaling functions and wavelets for $N=2$.}
We shall calculate values of coefficients of the obtained frequency
functions (19), (20) for various parameters $t$ and $u$ and we
shall find corresponding filters and wavelets $\varphi (x)$ and
$\psi(x)$. From formulas (21), (22) follows what enough to take
parameter values $t$ on interval of length $\pi$. We shall
consider the following parameter values
$t$: $0, \pm\pi/12, \pm\pi/6, \pm\pi/4, \pm\pi/3, \pm
5\pi/12, \pi/2$.

{\bf 4.1 Parameter values $t = 0$, $u = \pi/4$.} Coefficients of wavelets filters:
$$
h_0=\frac{1}{\sqrt{2}}(1,1,0,0), \qquad
g_1=\frac{1}{\sqrt{2}}(0,0,-1,1).
$$
We have obtained wavelets of Haar with the support on unit interval.
Refinement equations: $\varphi(x)=\varphi(2x)+\varphi(2x-1)$
and $\psi(x)=-\psi(2x-2)+\psi(2x-3)$.

{\bf 4.2. Parameter values $t = \pi/4$, $u =0$.} Coefficients of wavelets filters:
$$
h_0=\frac{1}{\sqrt{2}}(1,0,0,1), \qquad
g_1=\frac{1}{\sqrt{2}}(-1,0,0,1).
$$
We have obtained wavelets of Haar with the support on interval $[0,3]$.
Refinement equations: $\varphi(x)=\varphi(2x)+\varphi(2x-3)$
and $\psi(x)=-\psi(2x)+\psi(2x-3)$.

{\bf 4.3. Parameter values $t = \pi/2$, $u = –\pi/4$.}
Coefficients of wavelets filters:
$$
h_0=\frac{1}{\sqrt{2}}(0,0,1,1), \qquad
g_1=\frac{1}{\sqrt{2}}(-1,1,0,0).
$$
This is wavelets of Haar. Scaling function has the support
on interval $[2,3]$. Refinement equations:
$\varphi(x)=\varphi(2x-2)+\varphi(2x-3)$ and
$\psi(x)=-\psi(2x)+\psi(2x-1)$.

{\bf 4.4. Parameter values $t = –\pi/4$, $u = \pi/2$.}
Coefficients of wavelets filters:
$$
h_0=\frac{1}{\sqrt{2}}(0,1,1,0), \qquad
g_1=\frac{1}{\sqrt{2}}(0,1,-1,0).
$$
This is wavelets of Haar. Scaling function has the support
on interval $[1,2]$. Refinement equations:
$\varphi(x)=\varphi(2x-1)+\varphi(2x-2)$ and
$\psi(x)=\psi(2x-1)-\psi(2x-2)$.

{\bf 4.5. Parameter values $t=\pi/12$, $u =\pi/6$.}
Coefficients of wavelets filters:
$$
h_0=\frac{\sqrt{2}}{8}(3+\sqrt{3},1+\sqrt{3},1-\sqrt{3},3-\sqrt{3}),
\qquad
g_1=\frac{\sqrt{2}}{8}(-3+\sqrt{3},1-\sqrt{3},-1-\sqrt{3},3+\sqrt{3}).
$$
The result will be wavelets with coefficients which are obtained
by permutation of coefficients of the classical Daubechies
wavelets with the support of length 3. Refinement equations:
$$
\varphi(x)=\frac{3+\sqrt{3}}{4}\varphi(2x)
+\frac{1+\sqrt{3}}{4}\varphi(2x-1)
+\frac{1-\sqrt{3}}{4}\varphi(2x-2)
+\frac{3-\sqrt{3}}{4}\varphi(2x-3).
$$
In figure 1 graphs of wavelets are shown.

\begin{figure}
\epsfxsize154pt\epsfbox[0 0 180 165]{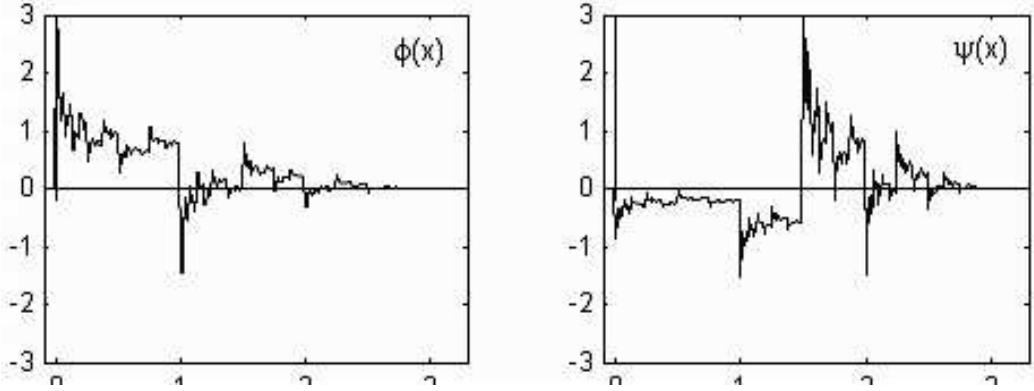}
\caption{Graphs of functions $\varphi(x)$ and $\psi(x)$ for $t=\pi/12$, $u =\pi/6$} \label{N2-1}
\end{figure}

{\bf 4.6. Parameter values $t=5\pi/12$, $u =-\pi/6$.}
Coefficients of wavelets filters:
$$
h_0=\frac{\sqrt{2}}{8}(3-\sqrt{3},1-\sqrt{3},1+\sqrt{3},3+\sqrt{3}),
\qquad
g_1=\frac{\sqrt{2}}{8}(-3-\sqrt{3},1+\sqrt{3},-1+\sqrt{3},3-\sqrt{3}).
$$
This example differs from previous only that coefficients of the
filter $\{h_n \}$ go upside-down. In this case scaling function
can be obtained from scaling function of example 4.5 with the
using of argument replacement: $\varphi(3-x)$. It follows from the
fact: if $\varphi (x)$ -- scaling function with the compact
support $[0, L]$ and the filter $\{h_n \}$ then function
$\varphi(L-x)$ also is scaling with the filter $\{h_{L-n} \}$.The
corresponding wavelet also can be obtained from previous as:
$-\psi (3-x)$. The graph of scaling function $\varphi(x)$ can be
obtained from the graph of the Fig.1 by mirroring about the line
$x=3/2$. For the graph of the wavelet $\psi(x)$ it is necessary to
add still mirroring about axis $Ox$ (Fig. 2).

\begin{figure}
\epsfxsize154pt\epsfbox[0 10 180 165]{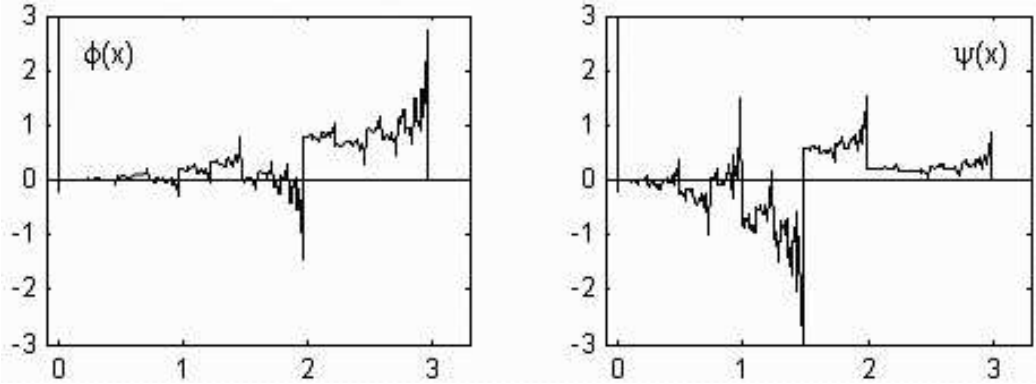}
\caption{Graphs of functions
$\varphi(x)$ and $\psi(x)$ for $t=5\pi/12$, $u =-\pi/6$ }
\label{N2-2}
\end{figure}

{\bf 4.7. Parameter values $t=-\pi/12$, $u = \pi/3$.}
Coefficients of wavelets filters:
$$
h_0=\frac{\sqrt{2}}{8}(1+\sqrt{3},3+\sqrt{3},3-\sqrt{3},1-\sqrt{3}),
\qquad
g_1=\frac{\sqrt{2}}{8}(-1+\sqrt{3},3-\sqrt{3},-3-\sqrt{3},1+\sqrt{3}).
$$
The result will be Daubechies wavelets with the support of length 3.
Refinement equation:
$$
\varphi(x)=\frac{1+\sqrt{3}}{4}\varphi(2x)
+\frac{3+\sqrt{3}}{4}\varphi(2x-1)
+\frac{3-\sqrt{3}}{4}\varphi(2x-2)
+\frac{1-\sqrt{3}}{4}\varphi(2x-3).
$$
In figure 3 graphs of wavelets are shown.

\begin{figure}
\epsfxsize154pt\epsfbox[0 0 180 165]{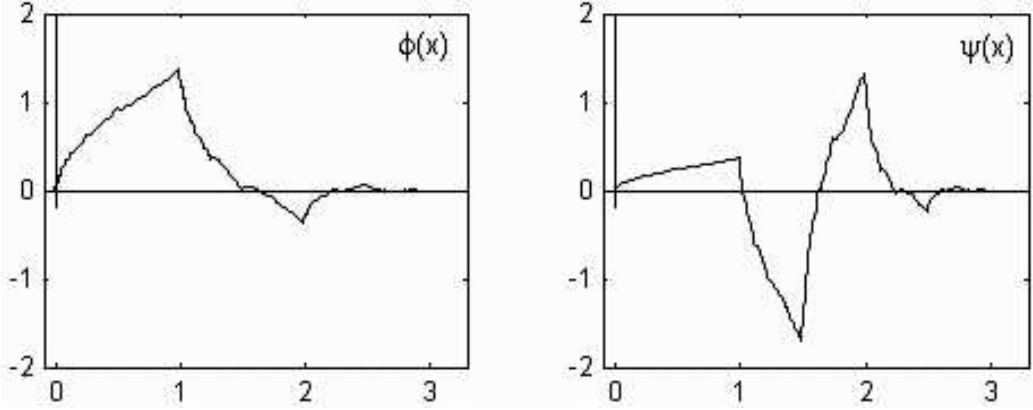} \caption{Graphs of functions
$\varphi(x)$ and $\psi(x)$ for $t=-\pi/12$, $u =\pi/3$}
\label{N2-3}
\end{figure}

{\bf 4.8. Parameter values $t = -5\pi/12$, $u = 2\pi/3$.}
Coefficients of wavelets filters:
$$
h_0=\frac{\sqrt{2}}{8}(1-\sqrt{3},3-\sqrt{3},3+\sqrt{3},1+\sqrt{3}),
\qquad
g_1=\frac{\sqrt{2}}{8}(-1-\sqrt{3},3+\sqrt{3},-3+\sqrt{3},1-\sqrt{3}).
$$
This example differs from the previous only that coefficients of
the filter $\{h_n \}$ go upside-down. In this case scaling
function can be obtained from Daubechies scaling function with the
help of argument replacement: $\varphi(3-x)$, and wavelet is
$-\psi(3-x)$.


{\bf 4.9. Parameter values $t = \pi/6$, $u = \pi/12$.}
Coefficients of wavelets filters:
$$
h_0=\frac{\sqrt{2}}{8}(3+\sqrt{3},3-\sqrt{3},1-\sqrt{3},1+\sqrt{3}),
\qquad
g_1=\frac{\sqrt{2}}{8}(-1-\sqrt{3},1-\sqrt{3},-3+\sqrt{3},3+\sqrt{3}).
$$
The result will be wavelets with coefficients which are obtained by
permutation of Daubechies wavelets coefficients.
The refinement equation:
$$
\varphi(x)=\frac{3+\sqrt{3}}{4}\varphi(2x)
+\frac{3-\sqrt{3}}{4}\varphi(2x-1)
+\frac{1-\sqrt{3}}{4}\varphi(2x-2)
+\frac{1+\sqrt{3}}{4}\varphi(2x-3).
$$
In figure 4 graphs of wavelets are shown.
\begin{figure}
\epsfxsize154pt\epsfbox[0 0 170 165]{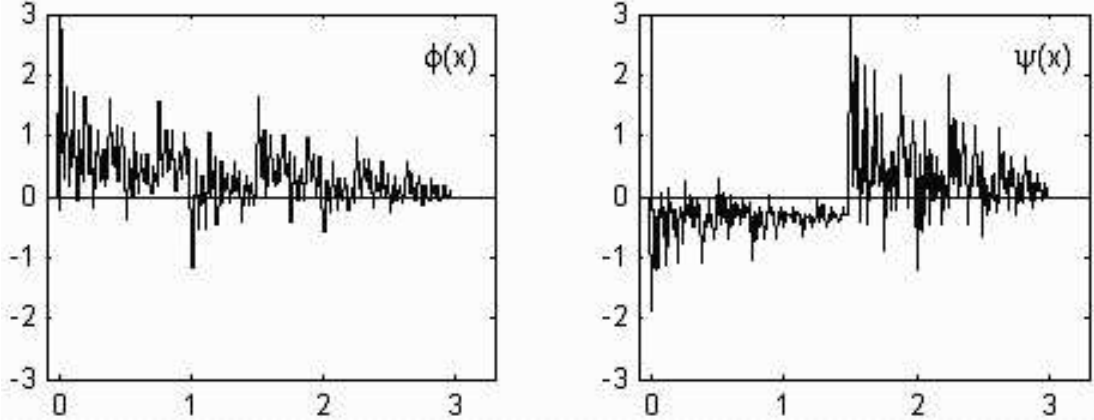} \caption{Graphs of functions
$\varphi(x)$ and $\psi(x)$ for $t = \pi/6$, $u = \pi/12$ }
\label{N2-4}
\end{figure}

{\bf 4.10. Parameter values $t = \pi/3$, $u = -\pi/12$.}
Coefficients of wavelets filters:
$$
h_0=\frac{\sqrt{2}}{8}(1+\sqrt{3},1-\sqrt{3},3-\sqrt{3},3+\sqrt{3}),
\qquad
g_1=\frac{\sqrt{2}}{8}(-3-\sqrt{3},3-\sqrt{3},-1+\sqrt{3},1+\sqrt{3}).
$$
This example differs from the previous only that coefficients of
the filter $\{h_n \}$ go upside-down. In this case scaling
function can be obtained from the previous scaling function by
replacement of argument: $\varphi(3-x)$, and wavelet is
$-\psi(3-x)$.

{\bf 4.11. Parameter values $t = -\pi/3$, $u = 7\pi/12$.}
Coefficients of wavelets filters:
$$
h_0=\frac{\sqrt{2}}{8}(1-\sqrt{3},1+\sqrt{3},3+\sqrt{3},3-\sqrt{3}),
\qquad
g_1=\frac{\sqrt{2}}{8}(-3+\sqrt{3},3+\sqrt{3},-1-\sqrt{3},1-\sqrt{3}).
$$
The result will be wavelets with coefficients which are obtained by
coefficients permutation of Daubechies wavelets with the support of
length 3. Refinement equations:
$$
\varphi(x)=\frac{1-\sqrt{3}}{4}\varphi(2x)
+\frac{1+\sqrt{3}}{4}\varphi(2x-1)
+\frac{3+\sqrt{3}}{4}\varphi(2x-2)
+\frac{3-\sqrt{3}}{4}\varphi(2x-3).
$$
In figure 5 graphs of wavelets are shown.
\begin{figure}
\epsfxsize154pt\epsfbox[0 0 180 125]{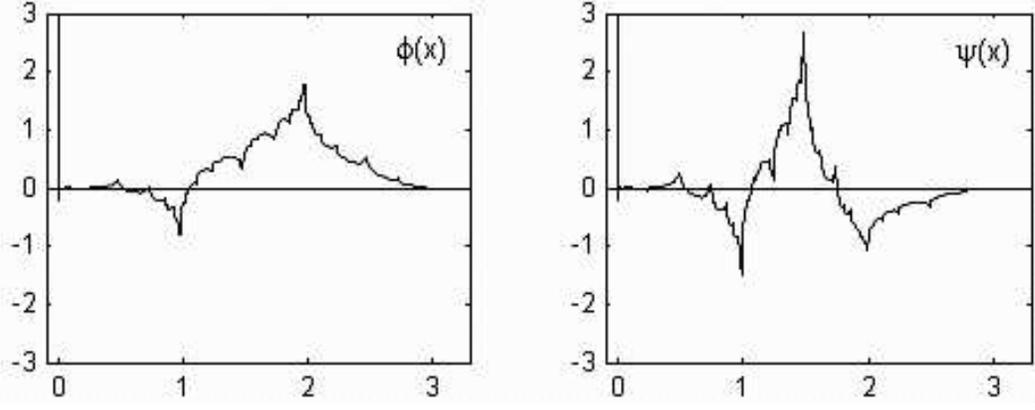} \caption{Graphs of functions
$\varphi(x)$ and $\psi(x)$ for $t = -\pi/3$, $u = 7\pi/12$ }
\label{N2-5}
\end{figure}

{\bf 4.12. Parameter values $t = -\pi/6$, $u = 5\pi/12$.}
Coefficients of wavelets filters:
$$
h_0=\frac{\sqrt{2}}{8}(3-\sqrt{3},3+\sqrt{3},1+\sqrt{3},1-\sqrt{3}),
\qquad
g_1=\frac{\sqrt{2}}{8}(-1+\sqrt{3},1+\sqrt{3},-3-\sqrt{3},3-\sqrt{3}).
$$
This example differs from the previous only that coefficients of
the filter $\{h_n \}$ go upside-down. In this case scaling
function and wavelet can be obtained from the previous by
replacement of argument: $\varphi(3-x)$, $-\psi(3-x)$.


{\bf 5. Construction of wavelets in case $N=3$.} In this section
we shall show the scheme of  scaling function and wavelets
construction for $N=3$. Though the diagonal matrix $D_k (w)$ can
be anyone, we shall take for example the diagonal matrix $D_1(w) =
{\rm diag} (1, w, 1)$, $|w |= 1$. The matrix $A_0 $ is:
$$
\left(
\begin{array}{ccc}
  1/\sqrt{3} & 1/\sqrt{3} & 1/\sqrt{3} \\
  1/\sqrt{2} & -1/\sqrt{2} & 0 \\
  1/\sqrt{6} & 1/\sqrt{6} & -2/\sqrt{6} \\
\end{array}
\right).
$$
Elements of the second orthogonal matrix $B_0$ should satisfy to conditions:
$$
b_{00}+b_{01}+b_{02}=1,\quad b_{10}+b_{11}+b_{12}=1,\quad
b_{20}+b_{21}+b_{22}=1.
$$
The solution of this system will be any set of orthonormal vectors
which coordinates satisfy to the equation of the plane
$x_0+x_1+x_2=1$. It is obvious, that coordinates of basis vectors
$e_1$, $e_2$, $e_3$ satisfy to this equation of plane. For this
solution the matrix $B_0$ it is identity. And we obtain the
wavelets of Haar,
$$
A_0 D_1(w)B_0= \left(
\begin{array}{ccc}
  1/\sqrt{3} & w/\sqrt{3} & 1/\sqrt{3} \\
  1/\sqrt{2} & -w/\sqrt{2} & 0 \\
  1/\sqrt{6} & w/\sqrt{6} & -2/\sqrt{6} \\
\end{array}
\right), \eqno(23)
$$
$$
H(z)= \frac{1}{\sqrt{3}}\left(
\begin{array}{ccc}
  1/\sqrt{3} & z^3/\sqrt{3} & 1/\sqrt{3} \\
  1/\sqrt{2} & -z^3/\sqrt{2} & 0 \\
  1/\sqrt{6} & z^3/\sqrt{6} & -2/\sqrt{6} \\
\end{array}
\right)
\left(
\begin{array}{ccc}
  1 & 1 & 1 \\
  z & \rho z & \rho^2 z \\
  z^2 & \rho^2 z^2 & \rho^4 z^2 \\
\end{array}
\right) ,
$$
$$
H_0(z)=\frac 13(1+z^2+z^4), \qquad H_1(z)=\frac
{1}{\sqrt{6}}(1-z^4), \qquad H_2(z)=\frac
{1}{3\sqrt{2}}(1-2z^2+z^4).
$$

The maximum degree of frequency function $H_0(z)$ is equal to
four, the support length $L$ is equal to two, as it is find from
the formula $L(N-1)+1 = \deg(H_0(z))+1$.

\begin{figure}
\epsfxsize154pt\epsfbox[0 0 180 165]{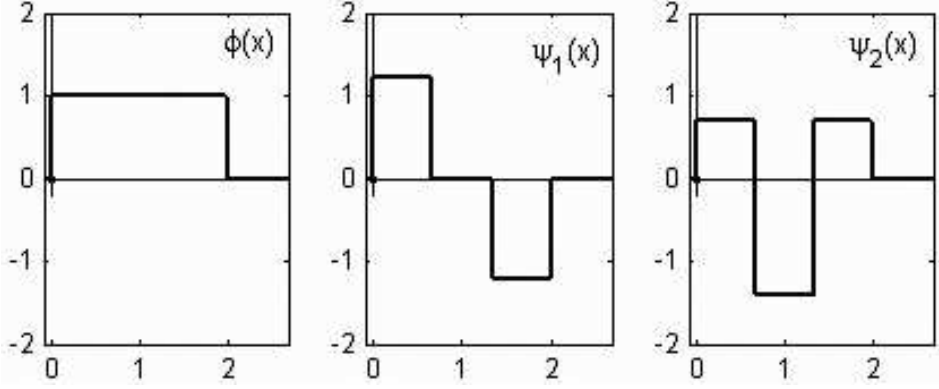} \caption{Graphs of functions
$\varphi(x)$, $\psi^1(x)$ and $\psi^2(x)$ for $t=0$} \label{pi0}
\end{figure}

It is easy to see, that scaling function $\varphi (x)$ is
characteristic function of interval [0,2), $\varphi (x) =
\chi_{[0,2)}(x)$. The refinement equation and wavelets (Fig.
\ref{pi0}):
$$
\varphi(x) = \varphi(3x)+\varphi(3x-2)+\varphi(3x-4),
$$
$$
\psi^1(x) = \frac{\sqrt{3}}{\sqrt{2}}\left(\varphi(3x) -
\varphi(3x-4)\right),
$$
$$
\psi^2(x) = \frac{1}{\sqrt{2}}\left(\varphi(3x) - 2\varphi(3x-2)
+\varphi(3x-4)\right),
$$

Any other solution can be obtained by rotation of basis vectors
$e_0$, $e_1$, $e_2$ in the plane $x_0+x_1+x_2=1$. We shall find
these solutions. Let
$$
M(t)=\left(
\begin{array}{ccc}
  1 & 0 & 0 \\
  0 & \cos t  & -\sin t \\
  0 & \sin t & \cos t \\
\end{array}
\right)
$$
-- the matrix of rotations around of the axis $e_0$. Then
$$
M_e(t)=A_0^{-1}M(t)A_0 = \frac 13 \left(\begin{array}{ccc}
  1+2\cos t & 1 -\cos t +\sqrt{3}\sin t  & 1 -\cos t -\sqrt{3}\sin t \\
  1 -\cos t -\sqrt{3}\sin t & 1+2\cos t  & 1 -\cos t +\sqrt{3}\sin t \\
  1 -\cos t +\sqrt{3}\sin t & 1 -\cos t -\sqrt{3}\sin t & 1+2\cos t \\
\end{array}
\right).
$$

Let's make rotation $M_e(t)e_k$ of column vectors $e_0 =(1, 0,
0)$, $e_1 =(0, 1, 0)$, $e_2 = (0, 0, 1)$, and we obtain rows of
the required matrix $B_0(t)$:
$$
B_0(t)= \frac 13 \left(
\begin{array}{ccc}
  1+2\cos t & 1 -\cos t -\sqrt{3}\sin t  & 1 -\cos t +\sqrt{3}\sin t \\
  1 -\cos t +\sqrt{3}\sin t & 1+2\cos t  & 1 -\cos t -\sqrt{3}\sin t \\
  1 -\cos t -\sqrt{3}\sin t & 1 -\cos t +\sqrt{3}\sin t & 1+2\cos t \\
\end{array}
\right). \eqno (24)
$$

Then $H (t, w) =A_0D_1 (w) B_0 (t) R (z)$
where the matrix $A_0 D_1 (z^N) $ is represented by the formula (23),
$B_0(t)$ -- by the formula (24) and the matrix $R (z)$ is
$$
R(z)=\frac{1}{\sqrt{3}}\left(
\begin{array}{ccc}
  1 & 1 & 1 \\
  z & \rho z & \rho^2 z \\
  z^2 & \rho^2 z^2 & \rho^4 z^2 \\
\end{array}
\right).
$$
Multiplying all these matrixes and choosing elements of the first column,
we obtain,
$$
H_0(t,z)=\frac 19 \left( 2+\cos t- \sqrt{3}\sin t +(2-2\cos t)z
+(2+\cos t +\sqrt{3}\sin t)z^2 + \qquad \quad \right.
$$
$$\left.
\qquad\qquad + (1-\cos t +\sqrt{3}\sin t)z^3 +(1+2\cos t)z^4
+(1-\cos t -\sqrt{3}\sin t)z^5 \right),\eqno (25)
$$
$$
H_1(t,z)=\frac {1}{3\sqrt{6}} \left( 1+2\cos t +(1-\cos t
-\sqrt{3}\sin t)z +(1-\cos t +\sqrt{3}\sin t)z^2 -  \qquad \right.
$$
$$\left.
\qquad\qquad -(1-\cos t +\sqrt{3}\sin t)z^3 -(1+2\cos t)z^4
+(-1+\cos t +\sqrt{3}\sin t)z^5 \right),\eqno (26)
$$
$$
H_2(t,z)=\frac {1}{9\sqrt{2}} \left(-1+4\cos t +2\sqrt{3}\sin t
-(1-\cos t +3\sqrt{3}\sin t)z -
 \qquad\qquad \qquad \   \right.
$$
$$
-(1+5\cos t +\sqrt{3}\sin t)z^2 +
$$
$$
\left.\qquad\qquad + (1-\cos t +\sqrt{3}\sin t)z^3 +(1+2\cos t)z^4
+(1-\cos t -\sqrt{3}\sin t)z^5 \right). \eqno (27)
$$

{\bf 6. Examples of scaling functions and wavelets for $N=3$.}
We shall calculate coefficients of the obtained frequency
functions (25), (26) and (27) for various parameter values $t$.
The obtained filters allow to find corresponding wavelets
$ \varphi (x)$, $\psi^1 (x)$ and $\psi^2 (x) $ by usual methods
\cite {Sm}, \cite {Db}. It is enough to find scaling function
$ \varphi (x)$. Wavelets - functions $\psi^1 (x)$ and $\psi^2 (x)$
are defined by formulas
$$
\psi^1(x)=\sqrt{N}\sum_{n\in \mathbb{Z}} g_n^1 \varphi(Nx-n),
\qquad \psi^2(x)=\sqrt{N}\sum_{n\in \mathbb{Z}} g_n^2
\varphi(Nx-n)
$$
with known filters $\{g_n^1 \}$ and $\{g_n^2 \}$ and function
$\varphi(x)$.

Let's consider the following parameter values $t$: 0,
$\pi/6$, $\pi/4$, $\pi/3$, $\pi/2 $, $2\pi/3$, $\pi $, $4\pi/3$.
For each case graphs of wavelets-functions are shown.


{\bf 6.1. Value of parameter $t=0 $.} This case has already been
considered above. It is wavelets of Haar with the support $[0,2]$
(Fig. \ref{pi0}).

{\bf 6.2. Value of parameter $t=\pi/6$.} Filters of scaling function
$\varphi (x)$ and wavelets $\psi^1 (x)$ and $\psi^2 (x)$:
$$
h_0=\frac{\sqrt{3}}{9}(2,2-\sqrt{3},2+\sqrt{3},1,1+\sqrt{3},1-\sqrt{3}),
$$
$$
g_1=\frac{\sqrt{6}}{18}(3+\sqrt{3},-3+\sqrt{3},\sqrt{3},-\sqrt{3},-3-\sqrt{3},3-\sqrt{3}),
$$
$$
g_2=\frac{\sqrt{6}}{18}(-1+3\sqrt{3},-1-\sqrt{3},-1-2\sqrt{3},1,1+\sqrt{3},1-\sqrt{3}).
$$
The refinement equation:
$$
\varphi(x) = \frac
13(2\varphi(3x)+(2-\sqrt{3})\varphi(3x-1)+(2+\sqrt{3})\varphi(3x-2)+\varphi(3x-3)+
$$
$$
+(1+\sqrt{3})\varphi(3x-4)+(1-\sqrt{3})\varphi(3x-5).
$$

\begin{figure}
\epsfxsize154pt\epsfbox[0 0 180 165]{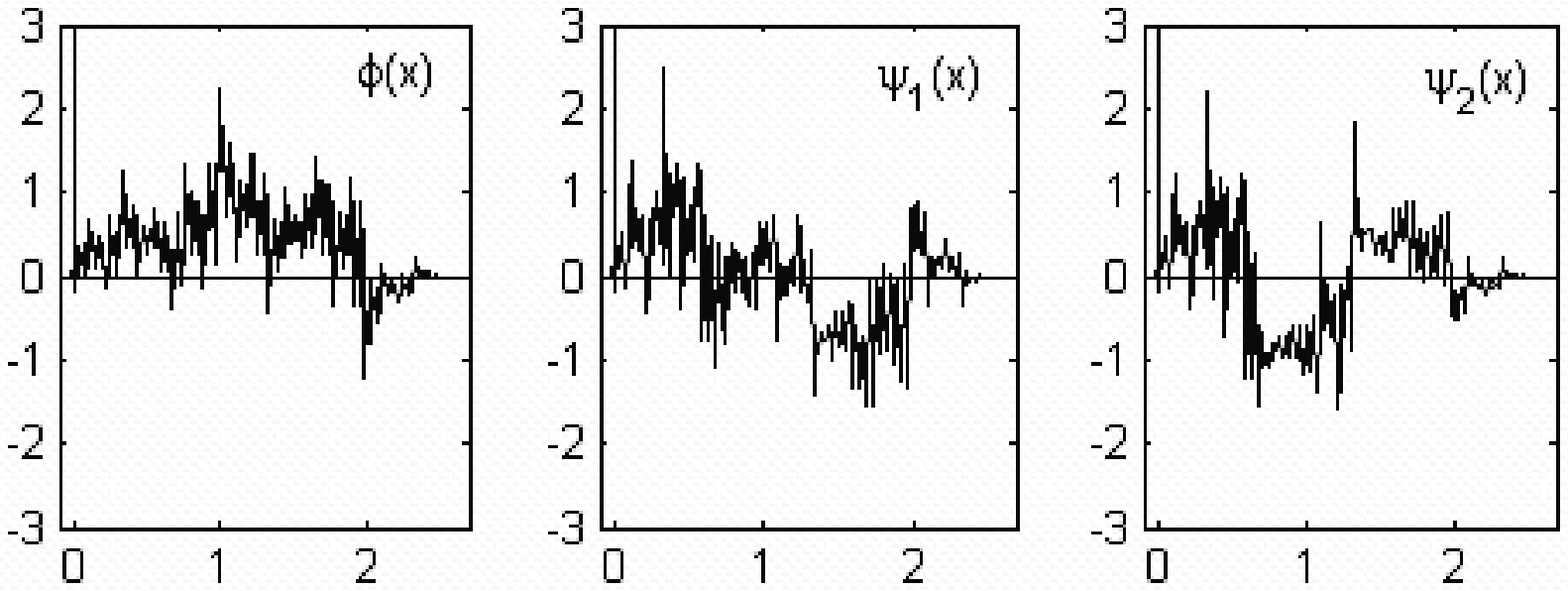}
\caption{Graphs of functions $\varphi(x)$, $\psi^1(x)$ and
$\psi^2(x)$ for $t=\pi/6$} \label{pi6}
\end{figure}
Graphs of wavelets are shown in figure \ref{pi6}

{\bf 6.3. Value of parameter $t=\pi/4$.} Filters of scaling function
$\varphi(x)$ and wavelets $\psi^1(x)$ and $\psi^2(x)$:
$$
h_0=\frac{\sqrt{3}}{18}(4+\sqrt{2}-\sqrt{6}, 4-2\sqrt{2},
4+\sqrt{2}+\sqrt{6}, 2-\sqrt{2}+\sqrt{6}, 2+2\sqrt{2},
2-\sqrt{2}-\sqrt{6}) ,
$$
$$
g_1=\frac{\sqrt{6}}{36}(2\sqrt{3}+2\sqrt{6},
-3\sqrt{2}+2\sqrt{3}-\sqrt{6}, 3\sqrt{2}+2\sqrt{3}-\sqrt{6},
-3\sqrt{2}-2\sqrt{3}+\sqrt{6},
$$
$$
\qquad\qquad\qquad\qquad\qquad\qquad -2\sqrt{3}-2\sqrt{6},
3\sqrt{2}-2\sqrt{3}+\sqrt{6}) ,
$$
$$
g_2=\frac{\sqrt{6}}{36}(-2+4\sqrt{2}+2\sqrt{6},
-2+\sqrt{2}-3\sqrt{6}, -2-5\sqrt{2}+\sqrt{6}, 2-\sqrt{2}+\sqrt{6},
$$
$$
\qquad\qquad\qquad\qquad\qquad\qquad 2+2\sqrt{2},
2-\sqrt{2}-\sqrt{6}).
$$

\begin{figure}
\epsfxsize154pt\epsfbox[0 0 180 165]{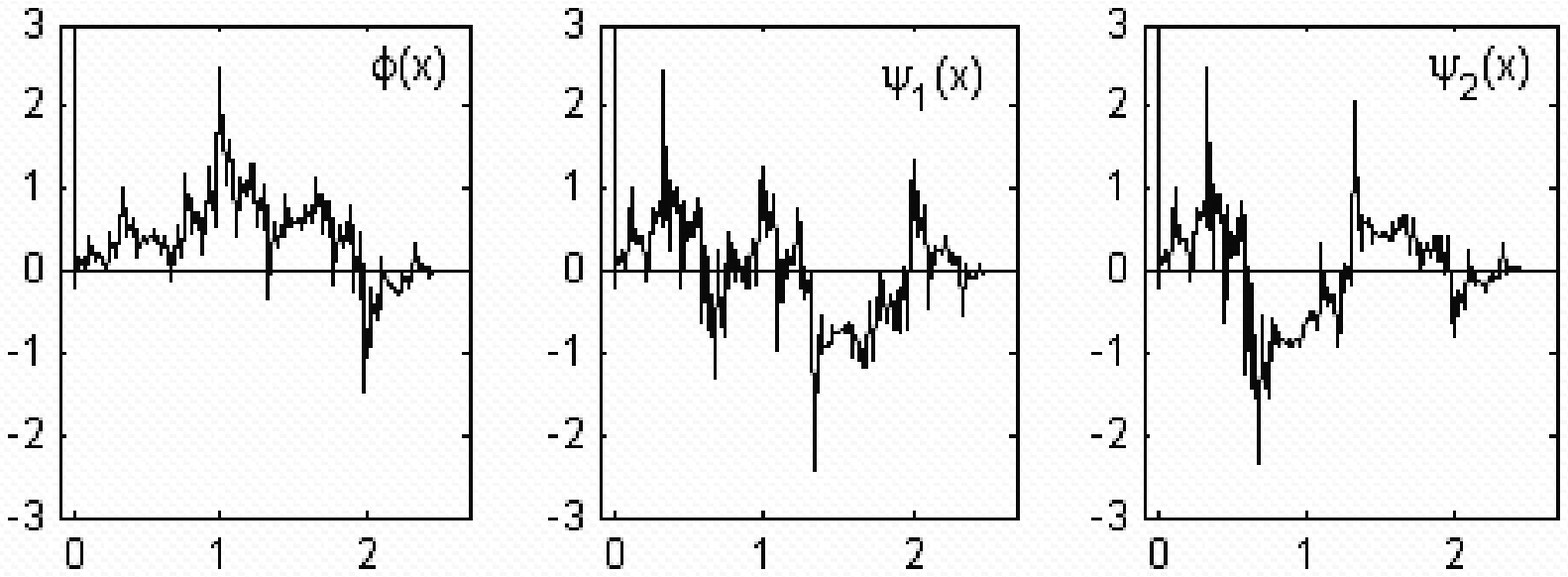}
\caption{Graphs of functions $\varphi(x)$, $\psi^1(x)$ and
$\psi^2(x)$ for $t=\pi/4$} \label{pi/4}
\end{figure}
Graphs of wavelets are shown in figure \ref{pi/4}

{\bf 6.4. Value of parameter $t=\pi/3$.} Filters of scaling function
$\varphi (x)$ and wavelets $\psi^1 (x)$ and  $\psi^2(x)$ are

$$
h_0=\frac{\sqrt{3}}{9}(1, 1, 4, 2, 2, -1) ,
$$
$$
g_1=\frac{\sqrt{2}}{6}(2, -1, 2, -2, -2, 1),  \qquad
g_2=\frac{\sqrt{6}}{18}(4, -5, -2, 2, 2, -1).
$$
The refinement equation:
$$
\varphi(x) = \frac
13(\varphi(3x)+\varphi(3x-1)+4\varphi(3x-2)+2\varphi(3x-3)+
2\varphi(3x-4)-\varphi(3x-5).
$$

\begin{figure}
\epsfxsize154pt\epsfbox[0 0 180 165]{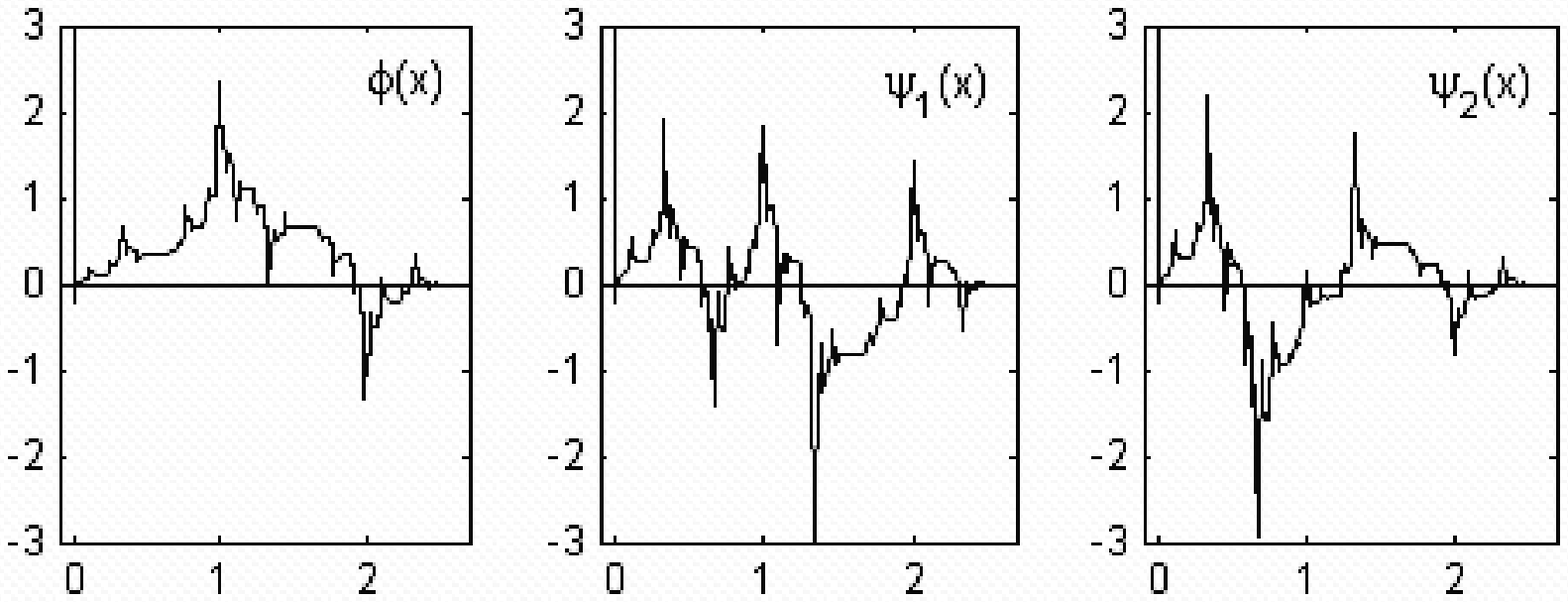}
\caption{Graphs of functions $\varphi(x)$, $\psi^1(x)$ and
$\psi^2(x)$ for $t=\pi/3$} \label{pi/3}
\end{figure}
Graphs of wavelets are shown in figure \ref{pi/3}.

{\bf 6.5. Value of parameter $t=\pi/2$.} Filters of scaling function
$\varphi (x)$ and wavelets $\psi^1 (x)$ and  $\psi^2(x)$:

$$
h_0=\frac{\sqrt{3}}{9}(2-\sqrt{3},2,2+\sqrt{3},1+\sqrt{3},1,1-\sqrt{3}),
$$
$$
g_1=\frac{\sqrt{6}}{18}(\sqrt{3},
-3+\sqrt{3},3+\sqrt{3},-3-\sqrt{3},-\sqrt{3},3-\sqrt{3}),
$$
$$
g_2=\frac{\sqrt{6}}{18}(-1+2\sqrt{3},-1-3\sqrt{3},-1+\sqrt{3},1+\sqrt{3},1,1-\sqrt{3}).
$$
The refinement equation:
$$
\varphi(x) = \frac
13((2-\sqrt{3})\varphi(3x)+2\varphi(3x-1)+(2+\sqrt{3})\varphi(3x-2)
+(1+\sqrt{3})\varphi(3x-3)+
$$
$$
+ \varphi(3x-4)+(1-\sqrt{3})\varphi(3x-5).
$$

\begin{figure}
\epsfxsize154pt\epsfbox[0 0 180 165]{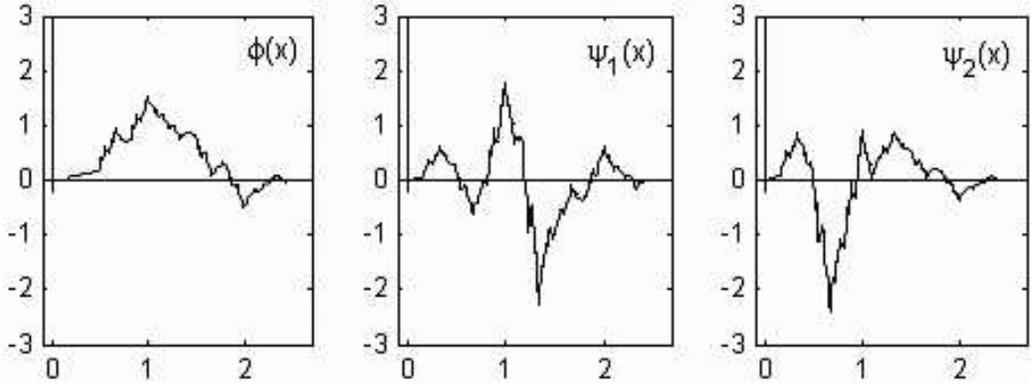}
\caption{Graphs of functions $\varphi(x)$, $\psi^1(x)$ and
$\psi^2(x)$ for $t=\pi/2$} \label{pi/2}
\end{figure}
Graphs of wavelets are shown in figure \ref{pi/2}.

{\bf 6.6. Value of parameter $t=2\pi/3$.} Filters of scaling function
$\varphi (x)$ and wavelets $\psi^1 (x)$ and  $\psi^2(x)$ is

$$
h_0=\frac{1}{\sqrt{3}}(0, 1, 1, 1, 0, 0) ,
$$
$$
g_1=\frac{1}{\sqrt{2}}(0, 0, 1, -1, 0, 0),  \qquad
g_2=\frac{1}{\sqrt{6}}(0, -2, 1, 1, 0, 0).
$$

It is wavelets of Haar. The scaling function
$\varphi(x)$ is characteristic function of interval $[1/2,3/2)$,
$\varphi(x)=\chi_{[1/2,3/2)}(x)$. The refinement equation and
wavelets:
$$
\varphi(x) = \varphi(3x-1)+\varphi(3x-2)+\varphi(3x-3),
$$
$$
\psi^1(x) =\frac{\sqrt{3}}{\sqrt{2}}
(\varphi(3x-2)-\varphi(3x-3)),
$$
$$
\psi^2(x) =
\frac{1}{\sqrt{2}}(-2\varphi(3x-1)+\varphi(3x-2)+\varphi(3x-3).
$$

\begin{figure}
\epsfxsize154pt\epsfbox[0 0 180 165]{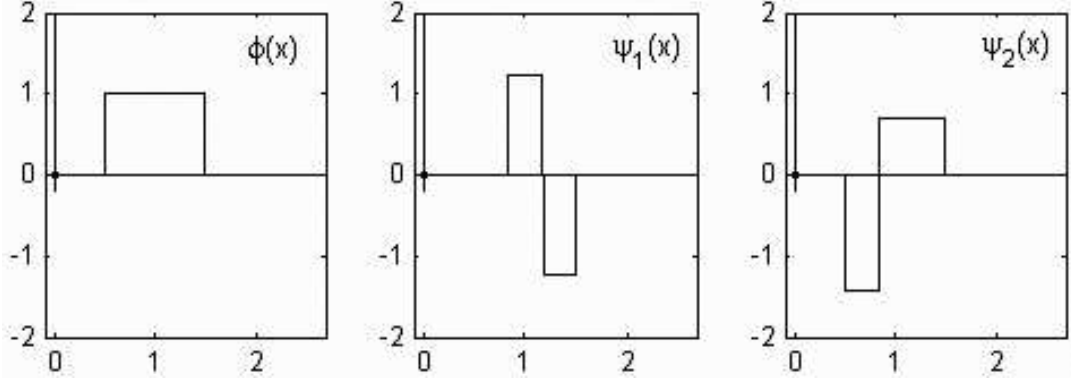}
\caption{Graphs of functions $\varphi(x)$, $\psi^1(x)$ and
$\psi^2(x)$ for $t=2\pi/3$} \label{2pi/3}
\end{figure}
In figure \ref{2pi/3} graphs of wavelets are shown.

{\bf 6.7. Value of parameter $t=\pi$.} Filters of scaling function
$\varphi (x)$ and wavelets $\psi^1 (x)$ and  $\psi^2(x)$:
$$
h_0= \frac{\sqrt{3}}{9}(1, 4, 1, 2, -1, 2) ,
$$
$$
g_1= \frac{\sqrt{2}}{6}(-1, 2, 2, -2, 1, -2),  \qquad g_2=
\frac{\sqrt{6}}{18}(-5, -2, 4, 2, -1, 2).
$$

The refinement equation:
$$
\varphi(x) = \frac 13 (\varphi(3x) +4\varphi(3x-1) +\varphi(3x-2)
+2\varphi(3x-3)-\varphi(3x-4) +2\varphi(3x-5).
$$

\begin{figure}
\epsfxsize154pt\epsfbox[0 0 180 165]{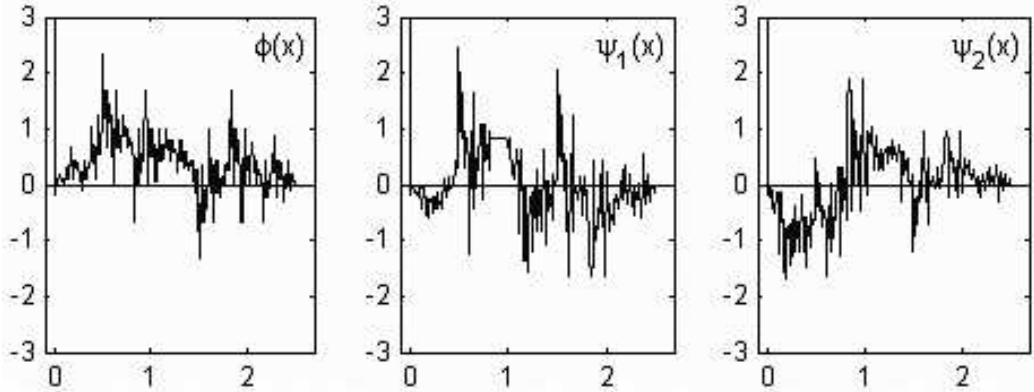}
\caption{Graphs of functions $\varphi(x)$, $\psi^1(x)$ and
$\psi^2(x)$ for $t=\pi$} \label{pi}
\end{figure}
Graphs of wavelets are shown in figure \ref{pi}.

{\bf 6.8. Value of parameter $t=4\pi/3$.} Filters of scaling function
$\varphi (x)$ and wavelets $\psi^1 (x)$ and  $\psi^2(x)$:
$$
h_0=\frac{1}{\sqrt{3}}(1, 1, 0, 0, 0, 1) ,
$$
$$
g_1=\frac{1}{\sqrt{2}}(0, 1, 0, 0, 0, -1), \qquad
g_2=\frac{1}{\sqrt{6}}(-2, 1, 0, 0, 0, 1).
$$

The refinement equation and frequency functions:
$$
\varphi(x) = \varphi(3x)+\varphi(3x-1)+\varphi(3x-5),
$$
$$
H_0(z)=\frac 13 (1+z+z^5),\qquad H_1(z)=\frac{1}{\sqrt{6}}
(z-z^5),\qquad H_2(z)=\frac{1}{3\sqrt{2}} (-2+z+z^5).
$$
Let's mark, that scaling function $\varphi(x)$ has a complicated
structure. Its support has fractal properties.

\begin{figure}
\epsfxsize154pt\epsfbox[0 0 180 165]{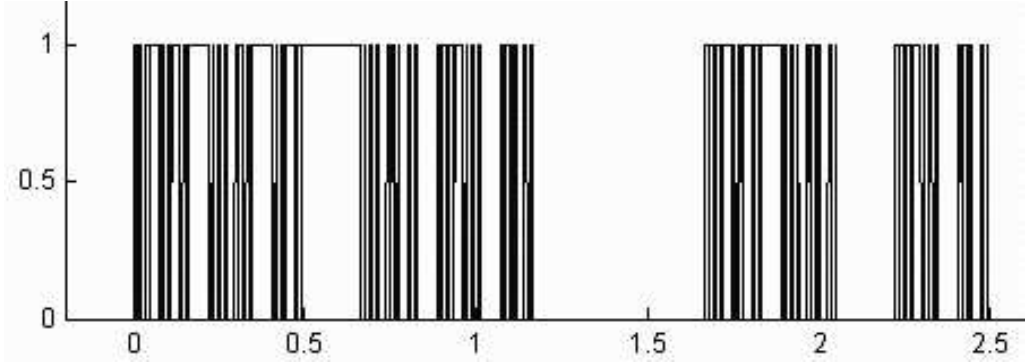}
\caption{Graph of function $\varphi(x)$ for $t=4\pi/3$} \label{4pi/3}
\end{figure}

\end{document}